\documentclass[a4paper, 11pt]{article}

\usepackage{graphicx}
\usepackage[all]{xy}
\usepackage{amsfonts, amssymb,url}
\usepackage{amsmath} 
\usepackage{amsthm}
\usepackage{epsfig}
\usepackage{ifthen}

\newtheorem{Thm}{Theorem}
\newtheorem{Lem}{Lemma}
\newtheorem{Prop}{Proposition}
\newtheorem{Coro}{Corollary}

\theoremstyle{definition}
\newtheorem{Rem}{Remark}
\newtheorem{Def}{Definition}
\newtheorem{Example}{Example}

\newcommand{\comment}[1]{}

\def\indd#1{{\bf 1}_{\{#1\}}}

\newcommand{\esp}{{\mathbb E}}

\newcommand{\defe}{\mathrel{\mathop:}=}

\newcommand{\inv}{^{-1}}

\newcommand{\var}{{\rm{Var}}}

\newcommand{\calB}{{\cal B}}

\newcommand{\filF}{{\cal F}}

\newcommand{\calN}{{\cal N}}

\def\indzd#1{\{#1_t\}_{t\in\mathbb Z^d}}

\newcommand{\eqnh}{\begin{eqnarray*}}
\newcommand{\eqne}{\end{eqnarray*}}
\newcommand{\eqnhn}{\begin{eqnarray}}
\newcommand{\eqnen}{\end{eqnarray}}
\newcommand{\equh}{\begin{equation}}
\newcommand{\eque}{\end{equation}}

\def\summ#1#2#3{\sum_{#1 = #2}^{#3}}
\def\prodd#1#2#3{\prod_{#1 = #2}^{#3}}

\def\sumzd#1{\sum_{#1\in\Zd}}

\newcommand{\widebar}{\overline}

\def\topp#1{^{(#1)}}

\def\nn#1{{\left\|#1\right\|}}

\def\snn#1{\|#1\|}

\def\bnn#1{\Big\|#1\Big\|}

\def\sabs#1{|#1|}
\def\babs#1{\Big|#1\Big|}
\def\ccbb#1{\left\{#1\right\}}
\def\bccbb#1{\Big\{#1\Big\}}

\def\bpp#1{\Big(#1\Big)}

\def\sbb#1{[#1]}
\def\bbb#1{\Big[#1\Big]}

\def\floor#1{\left\lfloor #1 \right\rfloor}
\def\sfloor#1{\lfloor #1 \rfloor}

\def\d{{\rm d}}
\def\e{{\rm e}}

\def\P{{\mathbb P}}
\def\B{{\mathbb B}}

\def\mand{\mbox{ and }}
\def\qmand{\quad\mbox{ and }\quad}

\def\qmwith{\quad\mbox{ with }\quad}
\def\mfa{\mbox{ for all }}

\def\mmas{\mbox{ as }}

\def\mfor{\mbox{ for }}

\def\adaptF#1{\{#1_t,\filF_t:0\leq t<\infty\}}

\def\weakto{\Rightarrow}

\def\limn{\lim_{n\to\infty}}
\def\limm{\lim_{m\to\infty}}

\def\liml{\lim_{l\to\infty}}

\def\wt#1{\widetilde{#1}}
\def\wb#1{\widebar{#1}}

\def\what#1{\widehat{#1}}

\def\Z{{\mathbb Z}}
\def\Zd{{\mathbb Z^d}}
\def\R{{\mathbb R}}
\def\Rd{{\mathbb R^d}}

\def\N{{\mathbb N}}

\def\ifhead#1#2{\left\{\begin{array}{#1@{\quad\mbox{ if }\quad}#2}}
\def\ifend{\end{array}\right.}

\def\indzd#1{\{#1_i\}_{i\in\Zd}}
\comment{

\def\P{\mathbb P}

\def\B{{\mathbb B}}

}

\begin{document}\sloppy

\title{An Invariance Principle for Fractional Brownian Sheets}


\author{Yizao Wang\vspace{0.4cm}\\
Department of Mathematical Sciences, University of Cincinnati}

\comment{
\institute{Y. Wang \at
              Department of Mathematical Sciences, University of Cincinnati,
2815 Commons Way, ML\#0025, Cincinnati, OH, 45221--0025, USA.\\              
              Tel.: +1--734--478--3642\\
              \email{yizao.wang@uc.edu}           
}

\date{Received: date / Accepted: date}
}

\maketitle

\begin{abstract}
We establish a central limit theorem for partial sums of stationary linear random fields with dependent innovations, and an invariance principle for anisotropic fractional Brownian sheets. Our result is a generalization of the invariance principle for fractional Brownian motions by Dedecker et al.~\cite{dedecker11invariance} to high dimensions. A key ingredient of their argument, the martingale approximation, is replaced by an $m$-approximation argument. An important tool of our approach is a moment inequality for stationary random fields recently established by El Machkouri et al.~\cite{elmachkouri13central}. 

\end{abstract}

\section{Introduction}
Consider a stationary linear random field $\indzd \xi$ with stationary mean-zero innovations $\indzd X$:
\equh\label{eq:model}
\xi_j = \sumzd i a_{j-i}X_i\,, j\in\Zd\,,
\eque
where
$\indzd a$ are a collection of real numbers such that $\sumzd ia_i^2<\infty$. 
In particular, we are interested in the partial sum over rectangles $\Lambda_n\subset\N^d$:
\[
S_n \defe \sum_{i\in\Lambda_n}\xi_i.
\]
We assume $\Lambda_n$ has the form $\Lambda_n = \{1,\dots,n_1\}\times\cdots\times\{1,\dots,n_d\}\subset\Zd$.
Set $b_{n,j} = \sum_{i\in\Lambda_n}a_{i-j}$ and $b_n = (\sumzd j b_{n,j}^2)^{1/2}$. Then $S_n = \sumzd jb_{n,j}X_j$. Throughout, we assume that the rectangle $\Lambda_n$ tends to $\N^d$ in the sense that for each $q=1,\dots,d$, $n_q \equiv n_q(n)\to\infty$ as $n\to\infty$.

We establish sufficient conditions for the following two problems:  
\begin{itemize}
\item[(i)] When do we have a {\it central limit theorem}
\[
\frac{S_n}{b_n} \equiv \frac{\sumzd jb_{n,j}X_j}{b_n}\weakto \calN(0,\sigma^2)\quad ?
\]
\item[(ii)] When do we have an {\it invariance principle}
\[
\bccbb{\frac{S_{\floor{nt}}}{b_n}}_{t\in[0,1]^d} \weakto \ccbb{\mathbb G_t}_{t\in[0,1]^d}
\]
and what is the limiting process $\mathbb G$?
\end{itemize}

These two problems have a long history. In the one-dimensional case ($d=1$), they  have been extensively investigated and many results are known. When $\{X_i\}_{i\in\Z}$ are independent and identically distributed (i.i.d.), see for example invariance principle results by Davydov~\cite{davydov70invariance} and Konstantopoulos and Sakhanenko~\cite{konstantopoulos04convergence}. When  $\{X_i\}_{i\in\Z}$ are stationary, see for example  Peligrad and Utev~\cite{peligrad97central,peligrad06central}, Wu and Woodroofe~\cite{wu04martingale}, and Merlev\`ede and Peligrad~\cite{merlevede06weak} for central limit theorems and Wu and Shao~\cite{wu06invariance}, and Dedecker et al.~\cite{dedecker11invariance} for invariance principles.
However, very few of the corresponding results in high-dimensional case ($d\geq 2$) are known, with the notable exceptions of Surgailis~\cite{surgailis82domains} and Lavancier~\cite{lavancier07invariance}. We compare our results to theirs in Section~\ref{sec:discussions}. We also point out in Remark~\ref{rem:lavancier} that the main results of \cite{lavancier07invariance} remain valid under a more general condition.

Recall that, for a given vector $H = (H_1,\dots,H_d)\in (0,1)^d$, a (real-valued) {\it fractional Brownian sheet} $\B^H = \{\B_t^H\}_{t\in[0,1]^d}$ with {\it Hurst index} $H$ is a real-valued mean-zero Gaussian random fields with covariance function given by
\equh\label{eq:covariance}
\esp(\B_s^H\B_t^H) = \prodd q1d \frac12\bpp{s_q^{2H_q} + t_q^{2H_q} - |s_q-t_q|^{2H_q}}, s,t\in[0,1]^d\,.
\eque
When $d = 1$, $\B^H$ is the fractional Brownian motion with Hurst index $H_1$.
 Fractional Brownian sheets play an important role in modeling anisotropic random fields with {\it long-range dependence} (see e.g.~Doukhan et al.~\cite{doukhan03theory} and Lavancier~\cite{lavancier06long}). They also arise in the study of stochastic partial differential equations (e.g.~Hu et al.~\cite{hu00stochastic} and~\O ksendal and Zhang~\cite{oksendal01multiparameter}). The investigation of their sample path properties is  another active research area (e.g.~Xiao~\cite{xiao09sample}).

This work develops an invariance principle for linear random fields with {\it stationary innovations}, converging to fractional Brownian sheets. Such a general consideration on innovations may have applications in spatial statistics (see e.g.~\cite{cressie93statistics}).
Our result can be seen as an extension of Dedecker et al.~\cite{dedecker11invariance} to high dimensions. The main difference is that we replace their martingale approximation method, which seems difficult to be generalized to high dimensions, by the $m$-approximation method. Besides, to characterize the weak dependence of innovations, we apply the {\it physical dependence measure} introduced by Wu~\cite{wu05nonlinear} and extended to random fields by El Machkouri et al.~\cite{elmachkouri13central}. In particular, El Machkouri et al.~\cite{elmachkouri13central} proved a moment inequality of weighted partial sums of $\indzd X$, which is of significant importance in the analysis of random fields by $m$-approximation. 

Another crucial assumption for our results is a product structure for the coefficients:
\[
a_i = \prodd q1d a_{i_q}(q)\,, i\in\Zd\,,
\]
and $\{a_{i_q}(q)\}_{i_q\in\Z}$ are square-summable real numbers for each $q=1,\dots,d$. In particular, the product structure allows us to extend the idea of {\it coefficient-averaging} by Peligrad and Utev~\cite{peligrad06central} to high dimensions. It also plays an important role in the analysis of asymptotic covariance structure.

The product structure of the coefficients is a reasonable assumption and it was assumed in~\cite{lavancier07invariance} to have an invariance principle for fractional Brownian sheets (the framework of \cite{lavancier07invariance} actually allows other forms of coefficients, although the limit may no longer be a fractional Brownian sheet). To see that it is not a restrictive assumption, recall the product form of the covariance formula~\eqref{eq:covariance}, and the fact that fractional Brownian sheets have the following stochastic integral representation with a product kernel $\prodd q1d g_{H_q}$:
\[
\B^H(t) = \frac1{\kappa_H}\int_{-\infty}^{t_1}\cdots\int_{-\infty}^{t_d}\prodd q1d g_{H_q}(t_q,s_q){\mathbb W}(\d s)\,,
\]
where $\mathbb W = \{\mathbb W_s\}_{s\in\Rd}$ is a standard Brownian sheet, for each $q$, $g_{H_q}(t_q,s_q) = ((t_q-s_q)_+)^{H_q-1/2} - ((-s_q)_+)^{H_q-1/2}$ with $r_+ = \max(r,0)$, and $\kappa_H$ is a normalizing constant (see e.g.~\cite{xiao09sample}).

The paper is organized as follows. We will provide the background on the physical dependence measure in Section~\ref{sec:prelim}. A central limit theorem (Theorem~\ref{thm:CLT}) and an invariance principle (Theorem~\ref{thm:fBs}) are established in Sections~\ref{sec:CLT}, and~\ref{sec:fBs} respectively. Discussions on related works are provided in Section~\ref{sec:discussions}.
\section{Preliminaries on physical dependence measure}\label{sec:prelim}
Consider stationary random fields $\indzd X$ of the following form
\equh\label{eq:Xi}
X_i = g(\epsilon_{i-j}:j\in\Zd), i\in\Zd\,,
\eque
where $g:\R^\Zd\to\R$ is a measurable function and $\indzd\epsilon$ are i.i.d.~random variables. Throughout this paper, we assume that $\esp X_0 = 0$.

El Machkouri et al.~\cite{elmachkouri13central} suggested to measure the dependence of $\indzd X$ as follows. Let $\epsilon^* = \indzd{\epsilon^*}$ be a random field coupled with $\epsilon$, defined by $\epsilon^*_i = \epsilon_i$ for all $i\in\Zd\setminus\{0\}$ and $\epsilon^*_0$ being a copy of $\epsilon_0$ independent of $\epsilon$. Set $X^*_i = g(\epsilon_{i-j}^*:j\in\Zd)$ and define the {\it physical dependence measure} of $\indzd X$ by
\equh\label{eq:Deltap}
\Delta_{p} \equiv\Delta_p(X) = \sum_{i\in\Zd}\nn{X_i-X_i^*}_p\,.
\eque
El Machkouri et al.~derived a central limit theorem and an invariance principle under the condition that $\Delta_p<\infty$ for certain $p\geq 2$. In particular, the following result is useful for our purpose.
\begin{Thm}[El Machkouri et al.~\cite{elmachkouri13central}]\label{thm:EMVW}
(i) Let $\indzd a$ be a family of real numbers. Then for any $p\geq 2$, 
\equh\label{eq:EMVW}
\bnn{\sum_{i\in\Zd}a_iX_i}_p\leq \bpp{2p\sum_{i\in\Zd}a_i^2}^{1/2}\Delta_p\,.
\eque
(ii) $\Delta_2<\infty$ implies that $\sumzd k|\esp X_0X_k|<\infty$.
\end{Thm}
\noindent In fact,~\eqref{eq:EMVW} was proved for $\indzd a$ having finite non-zero numbers, but the extension is immediate.

Our results require stationary random fields $\indzd X$ to satisfy $\Delta_p<\infty$ for some $p\geq 2$. El Machkouri et al.~\cite{elmachkouri13central} provided several such examples.
\begin{Example}
Consider $\indzd X$ in form of functional of linear random fields:
\[
X_i = g\bpp{\sumzd j\psi_{i-j}\epsilon_j}\,,i\in \Zd\,,
\]
where $g$ is a Lipschitz continuous function and the coefficients $\indzd\psi$ satisfy $\sumzd j|\psi_j|<\infty$. If $\epsilon_0\in L^p$ for $p\geq 2$, then $\Delta_p<\infty$ (\cite{elmachkouri13central}, Example 1). Note that this class of random fields include linear random fields.  Another class of non-linear random fields are the Volterra fields (\cite{elmachkouri13central}, Example 2).
\end{Example}
\comment{
\begin{Example}
Consider that in~\eqref{eq:Xi}, $X_i$ can be written as 
\[
X_i = \wt g(\epsilon_{i-j}: j\in\Gamma)
\] 
for some finite set $\Gamma\subset\Zd$. In this way, $\indzd X$ is an $m$-dependent stationary random field with some $m\in\N$. 
\end{Example}}

\begin{Rem}
Moment inequalities play an important role in establishing asymptotic results for random fields. Dedecker~\cite{dedecker01exponential} also established a similar moment inequality for random fields, under a different condition of weak dependence. 
In principle, in order to establish asymptotic normality one should expect to control, for finite subset $\Gamma\subset\Zd$, $\snn{\sum_{i\in\Gamma}X_i}_p\leq O(|\Gamma|^{1/2})$ for some $p\geq 2$.
Wang and Woodroofe~\cite{wang11new} established such an inequality for $\Gamma$ in form of rectangles, under a different condition of weak dependence.
\end{Rem}
\begin{Rem}
In the literature of (one-dimensional) stationary sequences, such a condition $\Delta_p<\infty$ on the weak dependence is often referred to as of {\it projective type}. An advantage of projective-type conditions is that, they often lead to easy-to-verify conditions for the asymptotic normality of various stationary processes arising from statistics and econometrics. For more on the projective-type conditions in one dimension, see for example Wu~\cite{wu05nonlinear,wu11asymptotic} and Merlev\`ede et al.~\cite{merlevede06recent}. For other types of conditions on weak dependence, see for example Bradley~\cite{bradley07introduction} and Dedecker et al.~\cite{dedecker07weak}.
\end{Rem}

\section{A central limit theorem}\label{sec:CLT}
We first establish a central limit theorem for triangular array in form of
\equh\label{eq:Sn}
S_n \equiv S_n(X,b)\defe \sumzd jb_{n,j}X_j\,,
\eque
with general {\it regular} coefficients $\{b_{n,j}\}_{n,j}$ to be defined below. For each $l\in\N$, set rectangle blocks of size $l^d$ by
\[
I_k \equiv I_k(l) = \{i\in\Zd: i_q \in\{lk_q+1,\dots,lk_q+l\}, q=1,\dots,d\}\subset \Zd\,,
\]
and define 
\equh\label{eq:cnk}
c_{n,k} = \frac1{l^d}\sum_{j\in I_k}b_{n,j}\,, n\in\N, k\in\Zd
\eque
and $c_n = (\sumzd kc_{n,k}^2)^{1/2}$. We introduce the following definition in the spirit of the coefficient-averaging idea of Peligrad and Utev~\cite{peligrad06central}.
\begin{Def}
We say the coefficients $\{\{b_{n,j}\}_{j\in\Zd}: n\in\N\}$ are {\it regular}, if $b_n^2\to\infty$ and for each $l\in\N$, 
\eqnhn
\limn\frac1{b_n^2}\sumzd k\sum_{j\in I_k}(b_{n,j}-c_{n,k})^2 & = & 0\label{eq:CS1}\\
\limn \frac1{b_n^2}\sumzd k\sum_{j\in I_k}|b_{n,j}^2-c_{n,k}^2| & = & 0\label{eq:CS2}\\
\limn\sup_{j\in\Zd} \frac{|c_{n,j}|}{c_n} & = & 0\label{eq:CS3}.
\eqnen
\end{Def}
Condition~\eqref{eq:CS2} implies $\limn l^dc_n^2/b_n^2 = 1$ which, combined with $b_n^2\to\infty$, yields $c_n>0$ eventually and thus~\eqref{eq:CS3} is well defined.
\begin{Thm}\label{thm:CLT}Consider $S_n$ as in~\eqref{eq:Sn} with some regular coefficients $\{b_{n,j}\}_{n,j}$. If $\Delta_2<\infty$, then, $\sigma^2\defe \sumzd k\esp X_0X_k<\infty$ and 
\equh\label{eq:CLT}
\frac{S_n}{b_n}\weakto\calN(0,\sigma^2)\,.
\eque
\end{Thm}
 Theorem~\ref{thm:CLT} leads to an answer to our first question. In particular, return to our problem of linear random fields with
$b_{n,j} = \sum_{i\in\Lambda_n}a_{i-j}$. It suffices to show such $\{b_{n,j}\}_{n,j}$ are regular. Recall that the coefficients $\indzd a$ are assumed to have the product structure:
\equh\label{eq:ai}
a_i = \prod_{q=1}^da_{i_q}(q), \mbox{ for some } \{a_{i_q}(q)\}_{i_q\in\Z}, q=1,\dots,d\,,
\eque
with $\sumzd ja_j^2<\infty$. This suffices to establish the regularity of $\{b_{n,j}\}_{n,j}$.
\begin{Coro}\label{coro:CLT}
If $b_{n,j} = \sum_{i\in\Lambda_n}a_{i-j}$ with $\indzd a$ having the product form~\eqref{eq:ai}, then $\{b_{n,j}\}_{n,j}$ is regular. As a consequence, if in addition $\Delta_2<\infty$, then $S_n(X,b)/b_n\weakto \calN(0,\sigma^2)$ with $\sigma^2 = \sumzd j\esp(X_0X_j)<\infty$.
\end{Coro}

We first prove Theorem~\ref{thm:CLT} through a series of approximations of $S_n$ in~\eqref{eq:Sn}. The main tool is Theorem~\ref{thm:EMVW}.\medskip

\noindent{\it (i): $m$-approximation.} First, fix $m$ and let $ \wb X = \indzd { \wb X}$ (depending on $m$) denote the stationary random field obtained by 
\[
\wb X_i = \esp(X_i\mid\filF^m_i)\,,
\]
where $\filF^m_i = \sigma(\epsilon_{i-j}: j\in\{-\floor{m/2},\dots,\floor{m/2}\}^d)$. The so-obtained random field $\indzd{\wb X}$ is $(m+1)$-dependent, that is, for all $i,j\in\Zd$, $\wb X_i$ and $\wb X_j$ are independent, if $\max_{q=1,\dots,d}|i_q-j_q|>m+1$. 

\begin{Lem}\label{lem:1} If $\Delta_p<\infty$ for some $p\geq 2$, then
\[
\lim_{m\to\infty}\sup_{n}\frac{\snn{S_n(X,b) - S_n(\wb X,b)}_p}{b_n}  = 0\,.
\]
\end{Lem}
\begin{proof}It follows from Proposition 3 and Lemma 2 in~\cite{elmachkouri13central}.
\end{proof}
In the sequel, let $\wb\Delta_p = \sumzd j\snn{\wb X_j - \wb X_j^*}_p$ denote the physical dependence measure of $\wb X$. Observe that $\Delta_p<\infty$ implies that $\wb\Delta_p = \sum_{j\in\{-m,m\}^d}\snn{\wb X_j-\wb X_j^*}_p <\infty$, for all $m\in\N$.\medskip

\noindent{\it (ii) Coefficient-averaging.}  This procedure was introduced by Peligrad and Utev~\cite{peligrad06central} in the one-dimensional case. 
For each $l\in\N$, 
recall the definition of $c_{n,k}$ in~\eqref{eq:cnk}. Set 
$\wb b_{n,j} \equiv \wb b_{n,j}(l) = \sumzd kc_{n,k}\indd{j\in I_k}$.

\begin{Lem}\label{lem:2} For each $m\in \N, l\in\N$, if $\wb\Delta_p<\infty$ for some $p\geq 2$, then
\[
\limn\frac{\snn{S_n(\wb X,b) - S_n(\wb X,\wb b)}_p}{b_n} = 0\,.
\]
\end{Lem}
\begin{proof}
Apply Theorem~\ref{thm:EMVW} and~\eqref{eq:CS1} to $\wb X$.
\end{proof}
\noindent{\it (iii) Big/small blockings.} Define, for $k\in\Zd, l\in\N, l>m+1$,
\[
\wt I_k \equiv \wt I_k(l) = \{i\in\Zd: i_q \in\{lk_q+1,\dots,lk_q+l-(m+1)\}, q=1,\dots,d\}\subset I_k\,.
\]
Set
$Y_k = \sum_{j\in \wt I_k} \wb X_k\,, k\in\Zd$. Since $\indzd {\wb X}$ are $(m+1)$-dependent, by the construction of $\{\wt I_k\}_{k\in\Zd}$, $\{Y_k\}_{k\in\Zd}$ are i.i.d.~random variables. Consider $S_n(Y,c) = \sumzd jc_{n,j}Y_j$. 
\begin{Lem}\label{lem:3}For each $m\in\N$, if $\wb\Delta_p<\infty$ for some $p\geq 2$ and~\eqref{eq:CS2} holds, then
\[
\liml\limn \frac{\snn{S_n(\wb X,\wb b) - S_n(Y,c)}_p}{b_n} = 0\,.
\]
\end{Lem}
\begin{proof}
Observe that $S_n(\wb X,\wb b) - S_n(Y,c) = \sumzd kc_{n,k}\sum_{j\in I_k\setminus\wt I_k}\wb X_j$. By Theorem~\ref{thm:EMVW}, 
\eqnh
\snn{S_n(\wb X,\wb b) - S_n(Y,c)}_p & \leq & \bccbb{2p{[l^d-(l-(m+1))^d]\sumzd kc_{n,k}^2}}^{1/2}\wb\Delta_p\\
& = & \bccbb{(2p)\bbb{1-\bpp{\frac{l-(m+1)}{l}}^d}(l^dc_n^2)}^{1/2}\wb\Delta_p\,.
\eqne
Note that~\eqref{eq:CS2} implies $\limn l^dc_n^2/b_n^2 = 1$, whence the desired result follows.
\end{proof}
\noindent{\it (iv) Triangular array of weighted i.i.d.~random variables.} Now we establish a central limit theorem for $S_n(Y,c)$. Recall that $Y$ depends on $l,m\in\N$.
\begin{Lem}\label{lem:4}For each $m\in\N, l>m+1$, if $\wb\Delta_p<\infty$ for some $p\geq 2$ and~\eqref{eq:CS3} holds, then
\[
\frac{S_n(Y,c)}{l^{d/2}c_n}\weakto \calN(0,\sigma_{m,l}^2)
\]
with 
\[
\sigma_{m,l}^2 = \sum_{i\in\{m+1-l,\dots,l-m-1\}^d}\esp(\wb X_0\wb X_i) \prod_{r=1}^d\bpp{1-\frac{m+1+|i_k|}l}\,.
\]
\end{Lem}
\begin{proof}
It is equivalent to prove a central limit theorem for $\sumzd j c_{n,j}\wb Y_j/c_n$ with $\wb Y_j = Y_j/l^{d/2}$. By straight-forward calculation, $\esp(\wb Y_1^2)/l^{d} = \sigma_{m,l}^2$. Then,~\eqref{eq:CS3} yields the desired result.
\end{proof}
\begin{proof}[Proof of Theorem~\ref{thm:CLT}]
Combining Lemmas~\ref{lem:1}--\ref{lem:4} yields that 
\[
\limm\lim_{l\to\infty}\limsup_{n\to\infty}\frac{\snn{S_n(X,b) - S_n(Y,c)}_2}{b_n} = 0\,.
\]
By~\cite{elmachkouri13central}, proof of Theorem 3.1 therein, $\sigma ^2 = \limm\lim_{l\to\infty}\sigma_{m,l}^2$. The desired result follows.
\end{proof}

At last, we prove Corollary~\ref{coro:CLT}.
By the product form~\eqref{eq:ai}, we can write
\equh\label{eq:bnj}
b_{n,j} = \prodd q1d b_{n,j}(q)\qmand b_n = \prodd q1d b_n(q)
\eque
with 
\equh\label{eq:bnjq}
b_{n,j}(q) \defe \summ i1{n_q} a_{i-j_q}(q)\qmand b_n(q) \defe \sbb{\sum_{j_q\in\Z} b_{n,j}^2(q)}^{1/2}.
\eque
 Observe that $b_{n,j}(q)$ only depends on $n_q, j_q\in\Z$ instead of $n,j\in\Zd$ (and we avoid writing $b_{n_q,j_q}(q)$ for the sake of simplicity). 
Accordingly, for fixed $l\in\N$ and $n,k\in\Zd$, we write $c_{n,k}(q) \defe \sum_{j_q=k_ql+1}^{(k_q+1)l}b_{n,j}(q)$ and $c_{n}(q) \defe [\sum_{k_q\in\Z} c^2_{n,k}(q)]^{1/2}$.

\begin{proof}[Proof of Corollary~\ref{coro:CLT}]
Fix $l\in\N$ and recall that $c_{n,k}$ depends on $l$. 
We first show $\sup_{j\in\Zd}|c_{n,j}|/c_n\to0$ as $n\to\infty$. Write
\equh\label{eq:supj}
\sup_{j\in\Zd}\frac{|c_{n,j}|}{c_n} \leq \prodd q1d\sup_{j\in\Zd}\frac{|c_{n,j}(q)|}{c_n(q)}\,.
\eque
It  suffices to show $\sup_{j\in\Zd}|c_{n,j}(q)|/c_n\to 0$ as $n\to\infty$ for all $q=1,\dots,d$. For each $q$, the convergence is implied by $\sup_{j\in\Zd}|b_{n,j}(q)|/b_n(q)\to 0$, proved by Peligrad and Utev~\cite{peligrad97central}, p.~448.

Next, we show that 
the product form~\eqref{eq:ai} implies that~\eqref{eq:CS1} and~\eqref{eq:CS2} hold. This result is an extension of Peligrad and Utev~\cite{peligrad06central}, Lemma A.1. 
We prove by induction. Note that now~\eqref{eq:CS2} becomes 
\equh\label{eq:CS2'}
\limn\frac1{b_n^2}\sumzd k\sum_{j\in I_k}\babs{\prodd q1d b_{n,j}^2(q)-\prodd q1dc_{n,k}^2(q)} = 0\,.
\eque
When $d=1$, it was shown in~\cite{peligrad06central} that~\eqref{eq:CS1} and~\eqref{eq:CS2'} holds for $\sum_{k\in\Z}a_k^2<\infty$. 

Suppose~\eqref{eq:CS1} and~\eqref{eq:CS2'} have been proved for $d-1$. We prove~\eqref{eq:CS1} for $d$ and the proof of~\eqref{eq:CS2'} is similar and omitted. By the inequality that for any real numbers $\alpha_q, \beta_q, q=1,\dots,d$, 
\[
\bpp{\prodd q1d\alpha_q - \prodd q1d\beta_q}^2 \leq 2\bbb{\bpp{\prodd q1{d-1}\alpha_q - \prodd q1{d-1}\beta_q}^2\alpha_d^2 + \bpp{\prodd q1{d-1}\beta_q^2}(\alpha_d-\beta_d)^2}\,,
\]
we bound $\sum_k\sum_{j}[\prodd q1db_{n,j}(q)-\prodd q1dc_{n,k}(q)]^2\leq 2(\Phi\topp1_n+\Phi\topp2_n)$ with
\[
\Phi_n\topp1 = {\sum_{k_1,\dots,k_{d-1}}\sum_{j_1,\dots,j_{d-1}}\bbb{\prodd q1{d-1}b_{n,j}(q) - \prodd q1{d-1}c_{n,k}(q)}^2}
 \times  \sum_{k_d}\sum_{j_d}b_{n,j}^2(d)
\]
and 
\eqnh
\Phi\topp2_n & = & \sum_{k_1,\dots,k_{d-1}}\sum_{j_1,\dots,j_{d-1}}\prodd q1{d-1}c^2_{n,k}(q) \sum_{k_d}\sum_{j_d}\bbb{b_{n,j}(d)-c_{n,k}(d)}^2\\
& = & l^{d-1}\prodd q1{d-1}c^2_n(q)\sum_{k_d}\sum_{j_d}\bbb{b_{n,j}(d)-c_{n,k}(d)}^2 \,.
\eqne
By induction, 
\[
\frac{\Phi_n\topp1}{b_n^2} = \frac{{\sum_{k\in\Z^{d-1}}\sum_{j\in\Z^{d-1}}\bbb{\prodd q1{d-1}b_{n,j}(q) - \prodd q1{d-1}c_{n,k}(q)}^2}}{\prodd q1{d-1}b_n^2(q)} = o(1)\,,
\]
and similarly, $\Phi_n\topp 2/b_n^2 \sim \sum_{k_d}\sum_{j_d}\sbb{b_{n,j}(d)-c_{n,k}(d)}^2/b_n^2(d) = o(1)$ by~\eqref{eq:CS1} with $d=1$. We have thus obtained~\eqref{eq:CS1} for all $d\in\N$. 
\end{proof}

\section{An invariance principle}\label{sec:fBs}
We consider weak convergence in the space $D[0,1]^d$ consisting of functions `continuous from above with limits from below' (see Bickel and Wichura~\cite{bickel71convergence} for details).
For $t\in[0,1]^d$, consider $S_n(t) \equiv \sum_{i_1=1}^{\floor{n_1t_1}}\cdots\sum_{i_d=1}^{\floor{n_dt_d}} \xi_i\in D[0,1]^d$.
This time we have
\equh\label{eq:Snt}
S_n(t) = \sumzd jb_{nt,j} X_j \qmwith b_{nt,j} \defe \sum_{i_1=1}^{\floor{n_1t_1}}\cdots\sum_{i_d=1}^{\floor{n_dt_d}}a_{i-j}. 
\eque
\begin{Thm}\label{thm:fBs}
Suppose there exists $H\in(0,1)^d$ such that
\equh\label{eq:bnq}
\limn \frac{b_{\floor {sn}}^2(q)}{b_n^2(q)} = s^{2H_q}, \mfa  s\in[0,1], q=1,\dots,d\,,
\eque
and there exists $p$ such that
\equh\label{eq:p}
p\geq 2,\ p > \max_{q=1,\dots,d} \frac1{H_q} \qmand \Delta_p<\infty\,.
\eque
Then, $\{S_n(t)/b_n\}_{t\in[0,1]^d}$ converges weakly in $D[0,1]^d$ to the fractional Brownian sheet with Hurst index $H$.
\end{Thm}
\begin{Rem}
When~\eqref{eq:bnq} holds for some $H\in(0,1)^d$ with $\max_{q=1,\dots,d}H_q\inv<2$, condition~\eqref{eq:p} becomes $\Delta_2<\infty$. Otherwise, we need to assume finite  higher-than-second-order moment to establish the tightness. A similar phenomena was observed in the one-dimensional case (\cite{dedecker11invariance}, Theorem 3.2). 
\end{Rem}
\begin{Example}\label{example:ai}
Due to the product structure, it suffices to provide examples of $\{a_{i_q}\}_{i_q\in\Z}$ such that~\eqref{eq:bnq} holds for each $q = 1,\dots,d$. Several examples have been provided in~Dedecker et al.~\cite{dedecker11invariance}, Examples 1--4. We summarize them below.
\begin{itemize}
\item[(i)] Fix $\alpha\in(0,1/2)$, and set $a_0 = 1, a_i = \Gamma(i+\alpha)/(\Gamma(\alpha)\Gamma(i+1))$ for $i\geq 1$. Then $H = \alpha+1/2$.
\item[(ii)] Fix $\alpha\in(0,1/2)$, and set $a_i = (i+1)^{-\alpha} - i^{-\alpha}$ for $i\geq 1$. Then $H = 1/2-\alpha$.
\item[(iii)] Fix $\alpha\in(1/2,1)$ and set $a_i\sim i^{-\alpha} l(i)$ for $i\geq 1$ with any slowly varying function $l$ at infinity. Then $H = 3/2-\alpha$.
\item [(iv)] Fix $\alpha>1/2$ and set $a_i\sim i^{-1/2}(\log i)^{-\alpha}$ for $i\geq 1$. Then $H = 1$.
\end{itemize}

\end{Example}
\begin{Example}[Fractionally integrated random fields]
Case (i) in Example~\ref{example:ai} above corresponds to the {\it fractionally integrated  random fields},  generated by {\it back-shift operators}. Let $B_q$ denote the back-shift operator on the $q$-th coordinate of the random fields $\indzd X$: $B_q X_i = X_{i_1,\dots,i_{q-1},i_q-1,i_{q+1},\dots,i_d}$. Then, for $\alpha_q\in(0,1/2)$,
\[
(I-B_q)^{-\alpha_q}X_j \defe \sum_{i=0}^\infty a_i(q)X_{j_1,\dots,j_{q-1},j_q-i,j_{q+1},\dots,j_d}\,, j\in\Zd\,,
\]
with $a_i(q)$ defined in Example~\ref{example:ai}, (i).
Thus,  fractionally integrated random fields defined by
\[
\xi_j = (I-B_1)^{-\alpha_1}\cdots(I-B_d)^{-\alpha_d}X_j\,, j\in\Zd
\]
fit in our model~\eqref{eq:model} with coefficients $a_i = \prodd q1d a_{i_q}(q), i\in \Zd$, where $\{a_{i_q}(q)\}_{i\in\Z_+}$ corresponds to $B_q$ as above and $a_{i_q}(q) = 0$ for $i_q<0$. This generalizes the fractional autoregressive integrated moving average (FARIMA) processes (see e.g.~\cite{wu06invariance} and references therein) to random fields (\cite{lavancier06long}). Note that Wu and Shao~\cite{wu06invariance} also established an invariance principle for the so-called {\it Type II} fractional integrated processes, which are slightly different from our model~\eqref{eq:model}. 
\end{Example}
Theorem~\ref{thm:fBs} follows as usual from the convergence of finite-dimensional distributions and tightness (see e.g.~Bickel and Wichura~\cite{bickel71convergence}), which are proved below separately.

\begin{Prop}[Convergence of finite-dimensional distributions]\label{prop:fdd}
Suppose $\Delta_2<\infty$ and for some $H\in(0,1)^d$~\eqref{eq:bnq} holds. 
then the finite-dimensional distributions of $\{S_n(t)/b_n\}_{t\in[0,1]^d}$ converge to that of a fractional Brownian sheet $\{\B_t^H\}_{t\in[0,1]^d}$ with Hurst index $H$.
\end{Prop}
\begin{proof}
Fix $m\in\N$. Take arbitrary $t\topp1,\dots,t\topp m\in [0,1]^d$ and write 
\[
n\topp r =(\sfloor{n_1t\topp r_1},\dots,\sfloor{n_dt\topp r_d})\in\Z_+^d \mfor r=1,\dots,m.
\]
Then we can write
\[
{\summ r1m\lambda_rS_n(t\topp r)} = {\sumzd j\bpp{\summ r1m{\lambda_rb_{n\topp r,j}}}\xi_j}\,.
\]
By Theorem~\ref{thm:CLT}, it suffices to show that $\{\wt b_{n,j} = \summ r1m \lambda_rb_{n\topp r,j}\}_{n,j}$ are regular and
\equh\label{eq:fdd}
\wt b_n^2\defe{\sumzd j\bpp{\summ r1m\lambda_rb_{n\topp r,j}}^2}\sim b_n^2\var\bbb{\summ r1m\lambda_r\B^H(t\topp r)}\,.
\eque

We first prove~\eqref{eq:fdd}. Observe that
\[
\wt b_n^2 = \sumzd j\summ r1m\summ s1m\lambda_r\lambda_sb_{n\topp r,j}b_{n\topp s,j}
= \summ r1m\summ s1m\lambda_r\lambda_s\sumzd jb_{n\topp r,j}b_{n\topp s,j} \,.
\]
Recall~\eqref{eq:bnj},~\eqref{eq:bnjq} and in particular that $b_{n\topp r,j}(q)$ involves only $n\topp r_q$ and $j_q$. Then,
\eqnh
\sumzd jb_{n\topp r,j}b_{n\topp s,j} & = & \sumzd j\prodd q1d b_{n\topp r,j}(q)b_{n\topp s,j}(q)
\\
& = & \sum_{j_1,\dots,j_{d-1}\in\Z}\prodd q1{d-1}b_{n\topp r,j}(q)b_{n\topp s,j}(q)\sum_{j_d\in\Z}b_{n\topp r,j}(d)b_{n\topp s,j}(d)\\
& = & \prodd q1d\sum_{j_q\in\Z}b_{n\topp r,j}(q)b_{n\topp s,j}(q)
\\
&= & \prodd q1d\sum_{j_q\in\Z}\frac12\bbb{b_{n\topp r,j}^2(q) + b_{n\topp s,j}^2(q) - (b_{n\topp r,j}(q)-b_{n\topp s,j}(q))^2}\\
& = & \prodd q1d\frac12\bbb{b_{n\topp r}^2(q) + b_{n\topp s}^2(q) - b_{|{n\topp r}-{n\topp s}|}^2(q)}\,.
\eqne
We have thus shown that
\[
\wt b_n^2 = \summ r1m\summ s1m\lambda_r\lambda_s\prodd q1d\frac12\bbb{b_{n\topp r}^2(q) + b_{n\topp s}^2(q) - b_{|{n\topp r}-{n\topp s}|}^2(q)}\,.
\]
On the other hand, by the covariance formula~\eqref{eq:covariance},
\begin{multline*}
\var\bbb{\summ r1m\lambda_r\B^H(t\topp r)} \\
= \summ r1m\summ s1m \lambda_r\lambda_s \prodd q1d \frac12\bbb{(t_q\topp r)^{2H_q} + (t_q\topp s)^{2H_q} - |t_q\topp r - t_q\topp s|^{2H_q}}.
\end{multline*}
Now~\eqref{eq:fdd} follows from~\eqref{eq:bnq} by recalling that $b_n^2 = \prodd q1d b_n^2(q)$.

Next we check that $\{\wt b_{n,j}\}_{j\in\Zd}$ are regular. Accordingly define $\wt c_{n,k} = \sum_{j\in I_k}\wt b_{n,j}/l^d$ and $\wt c_n = [\sumzd k\wt c_{n,k}^2]^{1/2}$. Observe that $\wt c_{n,k} =  \summ r1m \lambda_rc_{n\topp r,j}$. Then, conditions~\eqref{eq:CS1}, and~\eqref{eq:CS2} become
\equh\label{eq:wtbc}
\frac1{\wt b_n^2}\sumzd k\sum_{j\in I_k}(\wt b_{n,j}-\wt c_{n,k})^2 \to 0\mand \frac1{\wt b_n^2}\sumzd k\sum_{j\in I_k}\sabs{\wt b_{n,j}^2-\wt c_{n,k}^2} \to 0 \mmas n\to\infty\,.
\eque
The first part follows from the observation that $(\wt b_{n,j}-\wt c_{n,k})^2\leq m\summ r1m\lambda_r^2(b_{n\topp r,j}-c_{n\topp r,k})^2$ and  for each $r$, $\sumzd k\sum_{j\in I_k}(b_{n\topp r,j}-c_{n\topp r,k})^2/b_n^2\to 0$ as $n\to\infty$. To show the second part of~\eqref{eq:wtbc}, observe that 
\begin{multline*}
\sumzd k\sum_{j\in I_k}|\wt b_{n,j}^2-\wt c_{n,k}^2| \leq \sumzd k\sum_{j\in I_k}|\wt b_{n,j}- \wt c_{n,k}||\wt b_{n,j}+ \wt c_{n,k}|\\
 \leq  \bpp{\sumzd k\sum_{j\in I_k}|\wt b_{n,j}- \wt c_{n,k}|^2}^{1/2}\bpp{\sumzd k\sum_{j\in I_k}|\wt b_{n,j}+ \wt c_{n,k}|^2}^{1/2}\,,
\end{multline*}
where the first term in the last product is of order $o(\wt b_n)$ (by the first part of~\eqref{eq:wtbc}), while the second term $O(\wt b_n)$, whence~\eqref{eq:wtbc} follows.

At last, condition~\eqref{eq:CS3} becomes $\limn\sup_{j\in\Zd}{|\wt c_{n,j}|}/{\wt c_n} = 0$. To see this, observe that the second part of~\eqref{eq:wtbc} implies that $l^d\wt c_n^2\sim \wt b_n^2$. It then follows from~\eqref{eq:fdd} and~\eqref{eq:CS2} that 
$\wt c_n^2\sim C c_n^2$ for some constant $C>0$. The rest of the proof is similar to the control of~\eqref{eq:supj} and omitted. We have proved the regularity of $\{\wt b_{n,j}\}_{n,j}$ and thus the proposition.

\end{proof}
\begin{Prop}[Tightness]\label{prop:tightness}
If there exists $p$ such that~\eqref{eq:p} holds,
then the process $\{S_n(t)/b_n\}_{t\in[0,1]^d}$ is tight in $D[0,1]^d$.
\end{Prop}
\begin{proof}
We will apply Lavancier~\cite{lavancier05processus_TR}, Corollary 3. By slightly modifying the argument therein, it suffices to show that there exists constants $\beta>1, p>0, C>0$, such that for all $t\in(0,1)^d$ and $n$ large enough,
\equh\label{eq:tightness}
{\nn{S_n(t)}_p^p}\leq Cb_n^p \prodd q1d t_q^{\beta}\,.
\eque

By Theorem~\ref{thm:EMVW}, for $p\geq 2$,
\equh\label{eq:tightness1}
\nn{S_n(t)}_p = \bnn{\sumzd j b_{nt,j}X_j}_p\leq (2p)^{1/2}b_{nt}\Delta_{p} =(2p)^{1/2}\prodd q1d b_{nt}\topp q\Delta_{p} \,.
\eque
Observe that~\eqref{eq:bnq} means that $b_n(q)$ is regularly varying. Now, by Taqqu~\cite{taqqu79convergence}, Lemma 4.1, for all $\gamma_q>0$, there exists $C_q>0$ such that ${b_{nt}}/{b_n}\leq C_q t^{H_q-\gamma_q}$, uniformly on $[0,1]$ for $n$ larger than some $n_q$.
Now,~\eqref{eq:tightness1} can be controlled by, for $n$ large enough and some constant $C$,
\[
\nn{S_n(t)}_p^p\leq Cb_n^p\prodd q1d t^{p(H_q-\gamma_q)}\Delta_p\,.
\]
If $p$  satisfies~\eqref{eq:p}, then one can choose $\gamma_q>0$ small enough so that $p(H_q-\gamma_q)>1$. It then follows that~\eqref{eq:tightness} holds with $\beta = \min_q p(H_q-\gamma_q)>1$. We have thus proved the tightness.
\end{proof}

\section{Discussions}\label{sec:discussions}

We compare our results with Surgailis~\cite{surgailis82domains} and Lavancier~\cite{lavancier07invariance}. Both established general frameworks of weak convergence, including  invariance principles for fractional Brownian sheets as a special case. 

Surgailis~\cite{surgailis82domains} proved more general results in the sense that he considered general functionals of linear random fields. However, he only considered independent innovations. Furthermore, he assumed finite moments of any order of the innovations. At last, he considered only isotropic random fields with coefficients in form of $a_i = \ell(|i|)b(i/|i|)|i|^{-\beta},i\in\Zd\setminus\{0\}$, where $\ell$ is a slowly varying function and $b$ is a continuous function on the sphere.\medskip

Lavancier~\cite{lavancier07invariance} established invariance principles for linear random fields with stationary innovations. 
The weak dependence of innovations is characterized by the following assumption:\medskip

\noindent {\bf H1} The random field $\indzd X$ is centered and {\it weakly stationary} (i.e., with shift-invariant covariance structure) with bounded spectral density, and there exists a random field $\{\B_t\}_{t\in(0,\infty)^d}$ such that
\equh\label{eq:H1}
\bccbb{\frac1{n^{d/2}}\sum_{i_1=1}^{\floor {nt_1}}\cdots\sum_{i_d=1}^{\floor{nt_d}}X_i }_{t\in{(0,\infty)^d}}\stackrel{\rm f.d.d.}{\rightarrow} \ccbb{\B_t}_{t\in(0,\infty)^d}\,.\medskip
\eque
In addition, assumptions on the coefficients are given in terms of their Fourier transforms. 
Under these assumptions, results in~\cite{lavancier07invariance} are established through a spectral convergence theorem (in the spirit of \cite{lang00convergence}). 

Lavancier \cite{lavancier07invariance}'s results are more general in the sense that they do not assume the product structure of coefficients, and most importantly they cover various limiting objects, not necessarily in form of stochastic integrals.
In particular, the limiting objects in~\cite{lavancier07invariance} are described through a linear mapping 
\equh\label{eq:I}
I:L^2(\Rd)\to L^2(\Omega,\calB,\P),
\eque
constructed explicitly in the proof of Theorem 1 in \cite{lavancier07invariance}. 
As a special case, Theorem 5 in \cite{lavancier07invariance} established invariance principle for fractional Brownian sheets (see also Theorem 2 and Remark 5 therein). This special case of \cite{lavancier07invariance} is comparable to our results here.

We first compare our assumptions on the coefficients with the ones in~\cite{lavancier07invariance}, Theorem 5, where the product form of coefficients were also assumed.
For each $q$ and given the  coefficients $\{a_i(q)\}_{i\in\Z}$, consider $\what a\topp q\in L^2([-\pi,\pi])$ defined as 
\[
\what a\topp q(\omega) = \sum_{j\in\Z} a_j(q) \e^{-\sqrt{-1}j\omega}\,, q=1,\dots,d.
\]
Focus on $\what a\topp q(\omega_q)$ and omit the index $q$ from now on. 
In~\cite{lavancier07invariance}, Theorem 5 (see also Remark 6), it was assumed that
\equh\label{eq:equivalent}
\what a(\omega) \sim C|\omega|^{-\alpha}, \mbox{ as } \omega\to 0, \mbox{ for some } \alpha\in(0,1/2), C>0\,.
\eque
By results on trigonometric series (see e.g.~Zygmund~\cite{zygmund68trigonometric}, Chapter V, Theorems 2.6 and 2.24),~\eqref{eq:equivalent} is equivalent to $a_j\sim C_1j^{\alpha-1}$ as $j\to\infty$ for some constant $C_1>0$, which is a special case in Example~\ref{example:ai}, (iii) and this case covers the concrete example (i). The other cases of Example~\ref{example:ai} are not covered by~\eqref{eq:equivalent}. Thus, our assumptions on the coefficients are more general than~\cite{lavancier07invariance} in the case of invariance principles for fractional Brownian sheets. 
\medskip

At last, we compare assumptions on the weak dependence of innovations in two approaches. 
To establish an invariance principle for fractional Brownian sheets, \cite{lavancier07invariance} requires more than {\bf H1}. Namely, in addition it requires the linear mapping $I$ in~\eqref{eq:I} to be an isometry, so that the limiting objects can be interpreted as stochastic integrals. 
It was shown in \cite{lavancier07invariance}, Theorem 1 that the linear mapping $I$ is an isometry when  $\indzd X$ are strong white noise. However, this remains true under more general assumption on the innovations, as explained by the following remark.
\begin{Rem}\label{rem:lavancier}
The linear mapping $I$ in~\eqref{eq:I} is an isometry, when \medskip

\noindent {\bf H1'} The random field $\{X_k\}_{k\in\Zd}$ is centered and weakly stationary with bounded spectral density, and~\eqref{eq:H1} holds with $\{\B_t\}_{t\in(0,\infty)^d}$ being a standard Brownian sheet.\medskip

Under {\bf H1'}, the fact that $I$ is an isometry can be proved by almost the same proof as in \cite{lavancier07invariance} (this was proved in \cite{lavancier07invariance} when the innovations are the white noise). Indeed, it was proved in Theorem 1 there that $I$ is a well defined bounded operator. To show that $I$ is an isometry, it suffices to follow carefully the proof in \cite{lavancier07invariance}, first two lines in p.275, and replace $B_n(t)$ by a sequence of standard Brownian sheets $B_n^*(t)$. 
This generalization of \cite{lavancier07invariance} was discovered during several personal communications between the author and Fr\'ed\'eric Lavancier, after the current paper is finished. 
\end{Rem}

Our assumption on the weak dependence of innovations implies {\bf H1'}. Indeed, by Theorem~\ref{thm:EMVW}, (ii), $\Delta_2<\infty$ implies that the spectral density of $\indzd X$ is bounded. By~\cite{elmachkouri13central}, Proposition 4, $\Delta_2<\infty$ also implies~\eqref{eq:H1}. 

In conclusion of the comparison, our approach and Lavancier~\cite{lavancier07invariance}'s use completely different techniques; we have stronger assumptions on the weak dependence of the innovations, but weaker assumptions on the coefficients. It is interesting to investigate how large is the class of random fields that satisfy {\bf H1'} but not $\Delta_2<\infty$. \medskip


{\bf Acknowledgments} The author is grateful to Fr\'ed\'eric Lavancier for many inspiring discussions. The author is  grateful to Dalibor Voln\'y and Mohamed El Machkouri, for inviting him to visit Laboratoire de Math\'emathques Rapha\"el Salem,  Universit\'e de Rouen in June, 2011, and kindly showing him an early version of~\cite{elmachkouri13central}. The author would like to thank Stilian Stoev and Michael Woodroofe for helpful discussions. The author also thank anonymous referees for helpful comments.


\def\cprime{$'$} \def\polhk#1{\setbox0=\hbox{#1}{\ooalign{\hidewidth
  \lower1.5ex\hbox{`}\hidewidth\crcr\unhbox0}}}

\end{document}